\input amstex
\documentstyle{amsppt}

\loadeufb
\loadeusb
\loadeufm
\loadeurb
\loadeusm

\magnification =\magstep 1
\refstyle{A}
\NoRunningHeads

\topmatter
\title  On a conjecture of Shafarevich \endtitle
\author Robert  Treger \endauthor
\address Princeton, NJ 08540  \endaddress
\email roberttreger117{\@}gmail.com \endemail
\endtopmatter

\document  

\head
1.  Introduction
\endhead 
Unless stated otherwise,  $X \hookrightarrow  \bold P^{r}$ will be a nonsingular connected projective variety of dimension $n>0$. Let $U_X$ denote its universal covering.
Recall that the fundamental group $\pi_1(X)$ is {\it large}\/  if and only if $U_X$ contains no proper holomorphic subsets of positive dimension (Koll\'ar \cite{Kol}). 
 
\example{Conjecture (Shafarevich, 1972)} The universal covering of any nonsingular connected projective variety of  dimension $n>0$ with large fundamental group is a Stein manifold.
\endexample

The aim of this note is to prove the following
\proclaim{Theorem} If $\pi_1(X)$ is, in addition, residually finite then $U_X$ is a Stein manifold.
\endproclaim

In the sequel, we will assume that $\pi_1(X)$ is {\it residually finite}.

The conjecture attracted a great deal of attention, and recently it became a central problem in complex algebraic geometry. It is already known in several special cases (see Koll\'ar \cite{Kol}
%, a recent survey by Eyssidieux \cite{E} 
and references therein).
On the other hand,  Bogomolov and Katzarkov suggested   that the conjecture might fail in the case of nonresidually finite fundamental groups \cite{BK}.

 In Section 5, we will prove the theorem under an additional assumption that a general curvilinear section $C\subset X$ has the genus $g(C)\geq 2$ because, otherwise, $\pi_1(X)$ is  Abelian by the Campana-Deligne theorem \cite{Kol, Theorem\;2.14}, and the conjecture is well known when $\pi_1(X)$ is Abelian or, even, nilpotent  (Katzarkov \cite{Ka}). 

The idea of proof of the theorem is similar to the one by Siegel \cite{S}. He established  that
if  $U$ is a connected bounded domain in $\bold C^n$ covering a compact complex manifold $Y$ then $U$ is a domain of holomorphy. We will sketch his argument. 

He considers the Bergman metric on $U$  (see, e.g., \cite{Kob1}, \cite{Kob2, Chap.\;4.10}).  It is complete since $Y$ is compact. 
Recall the fundamental property of the Bergman metric, namely, it defines a natural  isometric embedding of $U$ into an infinite-dimensional projective space with a Fubini-Study metric. Infinite-dimensional projective spaces spaces were considered by Bochner \cite{B}, Calabi  \cite{C, Chap.\;4} and Kobayashi \cite{Kob1, Sect.\;7}. The idea of using square-integrable forms on arbitrary  manifolds can be found in Washnitzer \cite{W}.

 Let $\bold B(z,\bar z)$ denote the Bergman kernel of $U$.  Siegel proves that  $\log \bold B(z,\bar z)$
goes to infinity on any infinite discrete subset $T\subset U$. See also \cite{Kob1, Theorem 9.5}; Kobayashi observed in  \cite{Kob1, p.\;267} that his condition A.\!4 (\cite{Kob1, p.\;284} and \cite{Kob2, p.\;233, (C)}) is stronger than $\lim_{u\in T,u\to {\bar U}\backslash U} \log \bold B(z(u),\bar z(u))=\infty$. Hence  $U$ is a domain of holomorphy (equivalently, holomorphically compete or Stein domain)  by Oka's solution of the Levi problem. 

We observe that $Y$ is  a projective variety by the Poincar\'e ampleness theorem \cite{Kol, Theorem 5.22}. In 1950s, Bremermann proved  that an arbitrary bounded domain in $\bold C^n$ with complete Bergman metric is Stein \cite{Kob2, Theorem 4.10.21}. 

Now, Oka's solution of the Levi problem for domains in $\bold C^n$ admits a generalization due to Grauert (manifolds) and Narasimhan (complex spaces) (see \cite{N}). Thus, our aim is to define a metric on the manifold $U_X$ and a strictly plurisubharmonic function on $U_X$ that goes to infinity on an arbitrary infinite discrete subset of $U_X$.  The function generated by the diastasic potentials of the metric will be such a function.

\head
2.  Diastasis
\endhead 
The diastasis was introduced by Calabi \cite{C, Chap.\;2}.  
Let $M$ denote a complex manifold with real analytic Kahler metric. Let $\Phi$ denote a real analytic potential of the metric defined in a small neighborhood  $\Cal V \subset M$. Let $z=(z_1,\dots,z_n)$ be a coordinate system in $\Cal V$ and $\bar z=(\bar z_1,\dots,\bar z_n)$ a coordinate system in its conjugate neighborhood $\bar\Cal V\subset \bar M$. 
Let $(p,p)$ be a point on the diagonal of $M \times {\bar M}$ such that the neighborhood $\Cal V \times {\bar \Cal V} \subset M \times {\bar M}$ contains the point.

There exists a unique {\it holomorphic}\/ function $F$ on an open neighborhood of $(p,p)$ such that $F_{(p,p)}= \Phi_p$.
Here $\Phi_p$ is the germ at $p\in M$ of our real analytic function, and $F_{(p,p)}$ is the germ of the corresponding holomorphic function (complexification of $\Phi_p$ \cite {C, Chap.\;2}, \cite{U, Appendix}). One considers the sheaf $\Cal A^\bold R_M$ of germs of real analytic functions on $M$, and the sheaf $\Cal A^\bold C_{M \times \bar M}$ of germs of complex holomorphic functions on $M \times \bar M$. For each $p\in M$, we get a natural inclusion $\Cal A^\bold R_{M,p} \hookrightarrow \Cal A^\bold C_{M \times \bar M,(p,p)}$, called a complexification. The above equality is understood in this sense.

  Let $p$ and $q$ be two arbitrary points of $\Cal V$ with coordinates $z(p)$ and $z(q)$. Let $ F(z(p), \overline{z( q)})$ denote the corresponding complex holomorphic function on $\Cal V \times \bar \Cal V$. 
The {\it functional element of diastasis}\/ is defined as follows:
$$
D_M(p ,q):=F(z(p), \overline{z( p)})\,+\,F(z(q), \overline{z( q)})\,-\,F(z(p), \overline{z( q)})\,-\,F(z(q), \overline{z(p)}). 
$$
We get the germ $D_M(p,q)\in \Cal A^\bold C_{M \times \bar M,(p,q)}$, and $D_M(p,q)$ is uniquely determined  by the Kahler metric, symmetric in $p$ and $q$ and real valued  \cite{C, Prop.\;1,\,2}.
The {\it real}\/ analytic function generated by the above functional element is called the diastasis \cite{C, p.\;3}. 
The diastasis approximates the square of the geodesic distance in the {\it small}\/ \cite{C, p.\;4}.
For $\bold C^r$ with its unitary coordinates,
$
D_{\bold C^r}(p,q) = \sum^r_{i=1} |z_i(q) -z_i(p)|^2.
$

The fundamental property of the diastasis is that it is inductive on complex submanifolds \cite{C, Chap.\;2, Prop.\;6}.  
Another observation by Calabi \cite{C, {\it Note}\/ on p.\; 4} is that $D_M(p,q)$ represents the original functional element of the diastasis, and not values obtained from its maximal analytic continuation.

Now, let $\bold q\in M$ be a {\it fixed}\/ point, and $z=(z_1, \dots, z_n)$ a local coordinate system in a small neighborhood $\Cal V_\bold q \subset M$ with origin at $\bold q$. The real analytic function $\tilde \Phi_\bold q (z(p), \overline{z(p)}):= D_M(\bold q,p)$ on $\Cal V_\bold q$ is called the {\it  diastasic potential at}\/ $\bold q$ of the Kahler metric. It is strictly plurisubharmonic function in $p$ \cite{C, Chap.\;2, Prop.\;4}. 

The {\it prolongation}\/ over $M$ of  the germ of  diastasic potential $\tilde \Phi_\bold q (z(p), \overline{z(p)})$ at $\bold q$ is a {\it function}\/ $\bold P_M:=\bold P_{M,\bold q} \in  H^0(\Cal A^\bold R_M, M)$ such that,  for every $u\in M$,  $\bold P_{M}(u)$ coincides  with  $D_M(\bold q, u) $ meaning $D_M(\bold q, u) $, initially defined  in a neighborhood of $\bold q$, can be extended to the whole $M$. Moreover, the germ of $\bold P_M$ at $u$ is the  diastasic potential of our metric at $u$ (see (5.1) below). 

Clearly, the prolongation over $M$ is not always possible.  
Now, let  $\bold P^N$ be a projective space with the Fubini-Study metric. For $\bold q \in\bold P^N$, we consider Bochner canonical coordinates $z_1, \dots, z_N$ with origin at  $\bold q$ on the complement of a hyperplane at infinity.  By Calabi \cite{Chap.\;4, (27)}, 
$
	D_{\bold P^r}(\bold q,p)= \log\bigl(1+\sum^N_{\sigma=1}|z_\sigma(p)|^2 \bigl).
$
In the homogeneous coordinates $\xi_0,\dots,\xi_N  $, where $z_\sigma := {\xi_\sigma/\xi_0}$, we get
$$
D_{\bold P^N}(\bold q,p)= \log {\sum^N_{\sigma=0} |\xi_\sigma(p)|^2 \over |\xi_0(p)|^2}.
$$

On the other hand, let $U\subset \bold C ^n$ be a bounded domain.  Let $z_1, \dots, z_n$ be a local system of coordinates with origin at a point $\bold q\in U$.  By the characteristic property of the diastasic potential (vanishing of some partial derivatives; see  \cite{B, pp.\;180-181}, \cite{C, p.\;3,\;p.\;14}, and  \cite{U, Appendix} where this property is explicitly stated):
$$
\partial^{|I|}\bold P_{U, \bold q}(\bold q)/\partial z_I=\partial^{|I|}\bold P_{U,\bold q}(\bold q)/\partial \bar z_I=0 \quad (I:=\{i_1, \dots, i_n\}\; \text{where} \,\; i_1, \dots, i_n \geq 0),
$$
we get  $\log \bold B(z,\bar z)$ is the diastasic potential at $\bold q$ of the Bergman metric on $U$. Thus, $\bold P_{U,\bold q}=\log \bold B(z,\bar z)$ and it  is defined over the whole $U$, namely, $\bold P_{U,\bold q}(u)= D_U(\bold q, u)$.

\head
3. Bergman-type metrics in tower of coverings
\endhead

In this section, we will establish a version of a statement
attributed  to Kazhdan by Yau  \cite{Y, p.\;139}. For a survey of known results
and historical remarks, see a recent article by Ohsawa \cite{O, Sect.\;5}. In fact, Kazhdan had not presented a proof. In case $U_X$ is the disk $\Delta$, the first proof was given by Rhodes \cite{R}. Recently, McMullen has given a short proof for the disk \cite{M, Appendix}.
\smallskip

(3.1)  We consider a
tower of Galois coverings with each $Gal (X_i / X)$ a finite group:
$$
X=X_0 \leftarrow X_1 \leftarrow X_2 \leftarrow \cdots \leftarrow U, \quad \bigcap_i Gal(U/X_i) = \{1\} \; (0\leq i<\infty). \tag{3.1.1}
$$
We do not assume $U$ is simply connected. Let $\tau_i$ denote the projection $U\rightarrow X_i$, $\tau:=\tau_0$, and  $\tau_{jk}$ denote the projection $X_j\rightarrow X_k$ $(j\geq k)$.

The hyperplane bundle on $\bold P^{r}$ restricts to $\Cal L_X$,  called a polarization on  $X$.   Given our polarized Kahler metric  $g$ on  $X$, one can find a Hermitian metric  $h$ on $\Cal L_X$ with its Ricci curvature form equal to the corresponding Kahler form $\omega_g$. 
%\cite{Kob2, p.\;363}
%Recall the definition of Hermitian metric on $\Cal L$ \cite {Kol, Chap.\;7.1}. The %bundle $\Cal L\otimes \bar\Cal L$ has a natural conjugation $\delta: e\otimes e' %\rightarrow e'\otimes e$. Let $\eusb{R}\bold e(\Cal L\otimes \bar\Cal L)$ denote %the real subbundle fixed by $\delta$. The Hermitian metric on $\Cal L$ is a %vector bundle map $\Cal L\otimes \bar\Cal L\rightarrow \bold C_X$ which is real %%R \subset \bold C$. 
%In fact, $h= (\sum |\psi_\alpha|^2)$ where $\psi_\alpha$'s form a basis of %$H^0(X, \Cal L_X)$. 

We consider the volume form of the  Kahler form $\omega_g$. In local coordinates $z_1,\dots,z_n$ on $X$,
$
dv_g=V_g\! \prod^n_{\alpha=1}(\sqrt{-1}\!\cdot\!dz_\alpha\!\wedge\! d\bar {z}_\alpha)
$ 
where $V_g$ is a positive function. We will employ the same  volume form on all the coverings of $X$. Also, we will employ  the same Hermitian metric on all $\tau_{j0}^*(\Cal L_X)$.
\smallskip

(3.2) {\it Positive reproducing kernels and Bergman pseudometrics.}\/ The fundamental property of any  Berman-type pseudometric is the existence of a natural {\it continuous} map to a suitable projective space $\bold P(H^*)$ where the corresponding Hilbert space $H$ has a {\it reproducing kernel}.
\smallskip

(3.2.1) Let $M$
denote an arbitrary  complex manifold. Let $B(z,w)$ be a
Hermitian positive definite complex-valued  function  on
$M\times M$ which means:
\roster
 \item"{(i)}" $\overline{B(z, w)}=B(w,z), \qquad B(z,z)\geq 0;$
\item"{(ii)}" $\forall  z_1, \dots, z_N \in M, \quad
 \forall  a_1, \dots, a_N \in \bold C \implies 
\sum_{j,k}^N B(z_k,z_j)a_j \bar {a}_k \geq 0.$
\endroster
 If $B(z, w)$ is, in addition, holomorphic in the first variable then
$B$ is the reproducing kernel of a {\it unique Hilbert space}\/ $H$ of holomorphic functions on $U$ 
(see Aronszajn \cite{A, p.\;344,(4)} and the articles by Faraut and Kor\'anyi in \cite{FK, pp.\;5-14,
pp.\;187-191}). 
The evaluation at a point $Q\in M,$
$
e_Q:  f\mapsto f(Q),
$
is a continuous linear functional on $H$.
\smallskip

(3.2.2) Conversely, given a Hilbert space $H$ ($H\not =0$) of holomorphic functions
on
$M$ with all evaluation maps continuous linear functionals then, by the Riesz
representation theorem, for every $w\in M$ there exists a unique function $B_w\in H$
 such that $f(w)= \langle f, B_w\rangle$ ($\forall f\in H$) and $B(z,w):= B_w(z)$ is the reproducing kernel of
$H$ (which is Hermitian positive definite).

If we assume, {\it in addition}, that $B(z,z) > 0$ for every $z$ then  we can define   
$\log B(z, z)$ and a positive semidefinite  Hermitian form, called the Bergman pseudo-metric
$$
\qquad \qquad ds^2_M =2 \sum g_{j k} dz_j d\overline z_k, \qquad g_{j
k}:={\partial^2\log B(z, z)\over
\partial  z_j\partial\overline z_k  } .
$$
 We get a natural map $
\Upsilon : M \longrightarrow \bold {P}(H^*)$ whose image does not belong to a proper subspace of $\bold {P}(H^*)$   as in \cite{Kob2, Chap.\;4.10, pp.\;224-228}. As in  \cite{Kol, Chap.\;7, pp.\;81-84, Lemma-Definition 7.2}, the function $B(z,w)$ can be replaced  by a section of a relevant bundle. 
\smallskip

(3.3) Assuming $\eusm K_X$ is ample, we fix a large integer $t$ such that for every $i$,  $\eusm K^t_{X_i}$  is very ample (see \cite{Kol, 16.5}). The bundle $\eusm K^t_{X_i}$ is equipped with a Hermitian metric $h_{\eusm K^t_{X_i}}:=h^t_{\eusm K_{X_i}}$, where $h_{\eusm K_{X_i}}$ is the Hermitian metric  on $ \eusm K_{X_i}$ \cite {Kol, 5.12, 5.13, 7.1.1}. Also,  the bundle $\eusm K^t_{U}$ is equipped with a Hermitian metric.
%We get a standard volume form on $X$ corresponding to the very ample bundle %$\eusm K^t_X$.

Further, let $\psi_0,\dots ,\psi_N$ be an orthonormal basis of $H^0(X,\eusm K^t)$ with respect to $dv_g$ and $h_{\eusm K^t_X}$ (see also (4.3.1) below). Locally
$\psi_\beta =g_\beta  (z) (dz_1\wedge\cdots \wedge dz_n)^t$. We get an embedding $\sigma: X\hookrightarrow \bold P(H^0(X,K^t)^*)$.
We set 
$$
dv_{X,\eusm K^t}:=\big(\sum_{\beta =0}^N |g_\beta |^2 \big )^{1\over t}(\sqrt{-1})^{n^2} dz_1\wedge\cdots \wedge dz_n \wedge d\bar z_1 \wedge\cdots \wedge d\bar z_n.
$$
The associated Ricci form ${\text {Ric}}(dv_{X,\eusm K^t})$ is negative. If we pull back on $X$  the Fubini-Study metric on $\bold P(H^0(X,K^t)^*)$ then its Kahler form differs only by the sign from ${\text {Ric}}(dv_{X,\eusm K^t})$.

%We consider the same volume form on all the coverings of $X$.
\smallskip

(3.3.1) Let $\Omega_U^{(n,n)}$  denote the bundle of $(n,n)$-forms. Let $\Cal A \subsetneq U$  be a union of a countable number of real analytic subsets. As in \cite{Kol, Chap.\;7.1.1.2}, we fix a real homomorphism of $C^\infty$-bundles defined outside $\Cal A$:
$$
\Cal H_{U,t} : \eusm K_U^t\otimes \bar \eusm K_U^t \rightarrow \Omega_U^{(n,n)}\simeq\eusm K_U\otimes \bar \eusm K_U. 
$$
We consider the Hilbert space $H=H_{U}$ of  square-integrable
%, with  a weight $\rho >0$ invariant under the action of the Galois group, 
holomorphic weight $t$ differential forms $\omega$ on $U$. By square-integrable (or $L^2$),  we mean
$$
\int_U\Cal H_{U,t}(\omega\otimes \bar\omega)< \infty.
$$
 
We assume $H\not=0$. If all the evaluation maps are bounded (e.g., if $\Cal H$ is defined as in (3.3.2) below) then $H$ has a reproducing kernel as in the case of classical Bergman metric \cite{FK, pp.\;8-10, pp.\;187-188}.
Further, if the natural map
$$
 U\longrightarrow \bold P(H^*)
$$
is a holomorphic {\it embedding}\/ then the metric on $U$, induced from $\bold P(H^*)$, is called its $t$-{\it Bergman metric}. It is denoted by $b_{U,t}$ and the corresponding tensor is denoted by $g_{U,\eusm K^t}$. Of course, they depend on the choice of $\Cal H_{U,t}$.

Similarly, one defines the Euclidean space $V_i$ of  square-integrable 
holomorphic weight $t$  differential forms on $X_i$. 
 Let $\bold P(V^*_i) $ denote the corresponding projective space with its Fubini-Study metric. If the natural  map
$$
X_i \longrightarrow \bold P(V^*_i)
$$
is a holomorphic embedding then the induced metric on $X_i$ is called its $t$-{\it Bergman metric}. 
\smallskip

(3.3.2)   A sequence of homomorphisms $\Cal H_{X,t},\dots,\Cal H_{X_i,t},\dots\Cal H_{U,t}$ is said to be {\it compatible}\/ with the tower (3.1.1) if and only if $\tau_i^* \Cal H_{X_i,t}=\Cal H_{U,t}$ for all $i\geq 0$.

 In local coordinates, let $\omega_e=g_e(dz_1\wedge\cdots\wedge dz_n)^t$ be two weight $t$ forms ($e=1,2$).  Let $\varsigma$ be a suitable positive continuous function defined outside a  subset  $\Cal A \subsetneq X$, where $\Cal A$ is a union of a countable number of real analytic subsets. We assume  $\varsigma$ is bounded  away from $0$ on every compact subset of $X$. We consider inverse images of $\varsigma$ on all members of the tower. 
By abuse of notation, we denote the inverse images by the same symbol 
$\varsigma$. 
We can get a compatible sequence $\Cal H_{.,t}$ by setting 
$$
\Cal H_{.,t}(\omega_1,\omega_2):= (-2\sqrt{-1})^{-n}(-1)^{{n(n-1) \over 2}}\varsigma g_1 \bar g_2 dz_1\wedge\cdots \wedge dz_n\wedge d\bar z_1\wedge \cdots \wedge d\bar z_n.
$$ 

We consider $\Cal H_{.,t}$ defined as above with a suitable $\varsigma$.
The  corresponding Hilbert spaces have reproducing kernels as in the case of classical weighted Bergman spaces. Now, we assume, in addition, $\Cal L_X \subseteq\eusm K^t_X $.
Given the volume form $dv_g$ and the Hermitian metric $h$ on $\Cal L_X$, one can define an inner product on sections of $H^0(X,\Cal L)$ in two ways:
$$
\langle \omega_1, \omega_2\rangle' := \int_U h(\omega_1,\bar\omega_2)dv_g
\qquad  {\text {and}} \qquad
\langle \omega_1,\omega_2\rangle'' := \int_U \Cal H_{U,t}(\omega_1\otimes \bar\omega_2).
$$
As in \cite{Kol, Chap.\;5.1} with obvious modifications, one can compare the  corresponding norms. 
We have two Hilbert spaces with reproducing kernels, $H'$ and $H''\!.$ The  metrics on $U$, induced from $\bold P((H')^*)$ and $\bold P((H'')^*)$ respectively, will be equivalent.

\smallskip

(3.3.3) {\it Classical example.}\/ Let $U:=\Delta \subset \bold C$ be a disk. Let $\omega_1 :=g_1(z) (dz)^t$ and $\omega_2 :=g_2(z) (dz)^t$ be weight  $t$ holomorphic forms. With $\lambda_\Delta:=(1-|z|^2)^{-1}$, we  let  
$$
\Cal H_{\Delta,t} : \eusm K_\Delta^t\otimes \bar \eusm K_\Delta^t \longrightarrow \Omega_\Delta^{(1,1)}\simeq\eusm K_\Delta\otimes \bar \eusm K_\Delta, \qquad \omega_1\otimes \bar\omega_2 \mapsto \lambda_\Delta^{2-2t}g_1(z) \overline{g_2(z)}dz\wedge d\bar z. 
$$
As in (3.1.1), we consider a tower of Riemann surfaces.  In view of the correspondence between automorphic functions on $\Delta$ and differential forms on Riemann surfaces, we get the compatible sequence of  homomorphisms $\Cal H_{.,t}$.
\smallskip

The main result of this section is the following

 \proclaim{Proposition 1}  With the above notation, we assume that $U$ and all $X_i$'s have the $t$-Bergman metrics for an integer $t$, and the $\Cal H_{.,t}$'s,
defined with a help of the above $\varsigma$, are compatible with the tower. Then the
$t$-Bergman metric on $U$ equals the limit of pullbacks of the $t$-Bergman metrics from $X_i$'s.
\endproclaim

\demo{Proof} Let $b_{U,t}$ denote the
$t$-Bergman metric on $U$. Set $\tilde b_{U,t}:= \lim \sup  b_i $, where $b_i $ is the pullback on $U$ of the $t$-Bergman metric on $X_i$.  

  First, we get the inequality $\tilde b_{U,t} \leq b_{U,t}$.
We consider an open exhaustion of $U$, namely: $\{U_\nu \subset U,\; \nu=1, 2,\dots\}$, where each $\bar U_\nu$ is compact and $\bar U_\nu\subset U_{\nu+1}$.  Let $b(\nu)$ denote the $t$-Bergman metric on $U_\nu$. 
 It is well known that  $b_{U,t} = \lim_{\nu\to\infty} b(\nu)$. Indeed, a form is square-integrable on $U$ if and only if it is square-integrable on each $U_\nu$ and the integrals are bounded.
Given $U_\nu$, the restriction of $\tau_i$  on $U_\nu$ is one-to-one for $i\gg0$ because $\pi_1(X)$ is residually finite.
 Further,  $b(\nu) > b_i|U_\nu$ for all $i> i(\nu)$. This  establishes that $\tilde b_{U,t} \leq b_{U,t}$.

We have  $\tilde b_{U,t}=\lim b_i$. In the sequel, by a {\it local}\/ map we mean a map of a neighborhood. The above discussion suggests the following argument. 

 Let $V$ be the completion of the Euclidean space $E:=\cup V_i$. We will show that $U$, with the metric $\tilde b_{U,t}$, can be naturally isometrically embedded in $\bold P(V^*)$ with its Fubini-Study metric.

Let $p(w)\in U$ be a point and $\omega = f(dw_1\wedge\cdots \wedge dw_n)^t \in H$ a form
in local coordinates $w=(w_1, \dots, w_n)$ around $p(w)$. Recall that the $b_{U,t}$-isometric embedding $U\hookrightarrow \bold P(H^*)$ was given by the $t$-Bergman kernel section $B:= B_{U,\eusm K^t}$ (we do not assume  $\eusm K_U^t$ is trivial). The map $p(w) \mapsto \langle f, B_w\rangle$ has defined the local $b_{U,t}$-isometric holomorphic embedding
$U\rightarrow \bold P(H^*)$ (see \cite{Kob1, Sect.\; 7}, \cite{FK, pp.\;5-13 or p.\;188}) as well as the global $b_{U,t}$-isometric holomorphic embedding.

  Now, let $\omega = f(dw_1\wedge\cdots \wedge dw_n)^t $ be a form on $U$ that is the pullback of a form on $V_i$. As above, one can locally define  the map $p(w) \mapsto \langle f, B_w\rangle$ because $B=\lim_{\nu \to \infty} B(\nu)$ \cite{A, Part I, Sect.\;9}. The above map yields a  local  $b_i$-isometric holomorphic map $U\rightarrow \bold P(V_i^*)$. Hence we get a  local $\tilde b_{U,t}$-isometric holomorphic embedding  $U\hookrightarrow \bold P(V^*)$ (see \cite{FK, Prop.\; I.1.6 on p.\; 13}).

The local $\tilde b_{U,t}$-isometric holomorphic embedding yields the global one.
To show it, we use that $\pi_1(X)$ is residually finite. Indeed, the local embedding can be extended along a path \cite{C, Theorem 9}.  Since two paths with the same end points belong to a compact subset of $U$ that projects one-to-one to $X_i$ for $\forall i\gg 0$, we get the same result after extensions.

  Clearly, $U$ with the metric $b_{U,t}$ is $\bold P(H^*)$-resolvable at a point $p$ (1-resolvable of rank $N=\infty$ in Calabi's terminology \cite{C, p.\;19, Definition}), i.e,  the image of $U$ does not lie in a proper subspace of $\bold P(H^*)$.

Similarly, $U$ with the metric $\tilde b_{U,t}$ is $\bold P(V^*)$-resolvable at $p$. By the construction of $U\hookrightarrow \bold P(V^*)$, we obtain the natural embedding $\bold P(V^*)\hookrightarrow \bold P(H^*)$ induced by the natural map of the Hilbert spaces $H \rightarrow V$.  The latter embedding must be surjective hence $ \bold P(V^*) = \bold P(H^*)$.

This concludes the proof of the proposition.
\enddemo

\head
4. Metric $\Lambda$
\endhead

In this section, $\pi_1(X)$ is {\it not}\/ assumed to be large. We assume  $\pi_1(X)$ is residually finite and a general curvilinear section $C\subset X$ has $g(C) \geq 2$.
The metric $ \Lambda$ was suggested by a problem of Yau \cite{Y, Sect.\;6, p.\;139} who proposed to study $\lim_{t \to \infty}{1\over t} g_{X,\eusm K^t}$ when $\eusm K$ is the ample canonical bundle on $X.$ The metric $\Lambda$ is a generalization of the classical Poincar\'e metric though it is not necessary a Bergman metric. We will define a real analytic potential at every point of $U_X$.
\smallskip
(4.1)
First, we will consider the case: $C=X\hookrightarrow \bold P^r$ were $C$ is a connected nonsingular projective curve of genus $g(C)\geq 2$. We will assume the embedding is given by a very ample line bundle $\Cal L_C$ such that 
$$
\eusm K^{\ell}_C \subset \Cal L_C \subset \eusm K^{ m}_C,
$$
where $\eusm K_C$ is the canonical bundle and  $\ell, m$ are suitable  integers.
  We get Bergman-type metrics on $C$  (see, e.g., \cite{Y, Sect.\;6, p.\;138} and  \cite{Ti, p.\; 99}) and the Poincar\'e metric on $\Delta$. Since $\Delta$ is homogeneous,
$$
\bold B_{\Delta,\eusm K^t}(z,\zeta)= c(t)\bold B^t_{\Delta,\eusm K}(z,\zeta) \qquad  
(\bold B_{\Delta,\eusm K}(z,\zeta)=\pi^{-1}(1-z\bar\zeta)^{-2}),
\tag{4.1.1}
$$ 
where $t\gg 0$ is an integer, $c(t)$ is a known constant  depending on $t$ only, and $\bold B_{\Delta, \eusm K^t}$ denotes the $t$-Bergman kernel  (see, e.g., \cite{FK, p.\; 9}, \cite{Kol, (7.7.1)}). It follows
$$
\lim_{t \to \infty}{1\over t} g_{\Delta,\eusm K^t}=\biggl(\lim_{t \to \infty}{1\over t} {\partial^2\log\bold B_{\Delta,\eusm K^t}(z,z) \over \partial z \partial \bar z}\biggl)dz d\bar z =g_{\Delta,\eusm K}.\tag{4.1.2}
$$ 
\smallskip

(4.2) If $C$ is a general curvilinear section of  $X$ then we consider the inverse image of $C$ on $U_X$. By the Campana-Deligne theorem \cite{Kol, Theorem 2.14}, we obtain a connected open Riemann surface $R=R_C \subset U_X$ in place of the disk $\Delta$. We would like to construct a metric on $U_X$ such that its restriction on $R$ is well understood.

Set $\Gamma:= Gal(\Delta /R)$. Let 
$$
\bar \Cal F:= \big\{z\in \Delta \big|\, |Jac_\gamma(z)| \leq 1, \gamma\in \Gamma \big\}, \qquad \big\{\Cal F:= \{z\in \bar\Cal F \big|\, |Jac_\gamma(z)| < 1, \forall\gamma \not = 1\big\}
$$ 
be the fundamental domain of $R$ and the interior of the fundamental domain. 
One can also employ the Dirichlet-type fundamental domain centered at $\bold q:=0 \in \Delta$:
$$
\Cal D_\bold q(\Gamma):=\big\{ z\in \Delta \big|\; D_\Delta(z,\bold q)\leq D_\Delta(z,\gamma(\bold q)),\; \forall\gamma\in \Gamma\}.
$$
Indeed, $ D_\Delta(q_1,q_2) =  D_\Delta(\gamma q_1,\gamma q_2)$ for $\forall \gamma\in \Gamma$ because of the natural embedding of $\Delta$ into an infinite-dimensional projective space where  $\Gamma$ acts by collineations.
\smallskip
 
(4.2.1) High powers of $\Cal L_C:=\Cal L_X |_C$ are squeezed between powers of the canonical bundle on $C$.   
 For $t\geq1$, let $b_{R,t}$ denote the $t$-Bergman metric  on $R_C$ with $\Cal H_{.,t}$ as in the classical example (3.3.3).

%Here we consider weight $t$ differential forms. 
%The volume form $dv:=dv_g$ comes from the embedding of $R_C$ into infinite-
%dimensional projective space with the Fubini-Study metric \cite{R}. In fact, $dv$ is the %volume form of the classical Bergman metric on $R_C$. 
%The later metric is equal to $b_{R,1}$.

Let $D_{b,R,t}$ denote the functional element of  diastasis of $b_{R,t}$ at an arbitrary point of $R_C$.
We lift $b_{R,t}$ and $D_{b,R,t}$ to $\Delta$ and  get
$D_{b,R,t} \leq D_{b,\Delta,t}$ (locally at an arbitrary point of $\Delta$) by Proposition 1. Furthermore, it follows the convergence of the corresponding  holomorphic functions on $R_C\times \bar R_C$.
Hence
 $\lim_{t\to\infty}{1\over t} D_{b,R,t} \leq  D_{b,\Delta,1}.$ Set $b_{R}:=\lim _{t \to \infty}{1\over t}b_{R,t}.$ 

 It will be a  real analytic $Gal(R/C)$-invariant Bergman metric on $R_C$. Indeed, $b_{R,1}$ and ${1\over k}b_{R,t}$ are real analytic Kahler metrics. According to Tian \cite{Ti, Sect.\;4}, ${1\over t}b_{R,t}$ converge to $b_{R,1}$ hence  $b_{R}=b_{R,1}$. We get the metric $b_{R}$ whose diastasic potential $\bold P_{b,R}:=\bold P_{b,R,\bold a}$ is a {\it global function} on $R_C$, where $\bold a \in R_C $ is the image of the origin $0\in \Delta$.

One can replace $b_{R,1}$ by $b_{R,m}$, where $m$ is a sufficiently large {\it fixed}\/ number, and repeat the previous argument with $b_{R,m,t}$ in place of $b_{R,t}$. As before, we obtain 
$$
b_{R,m}= \lim _{t \to \infty}{1\over t}b_{R,m,t}.
$$ 
We will denote the functional element of the diastasis of $b_{R,m,t}$ by $D_{b,R,m,t}$.
\smallskip

(4.2.2) Now, let $\Cal L _R$ denote the inverse image of $\Cal L_C$ on $R_C$. We consider the metrics $g_{R,t}$ constructed, as in \cite{Ti}, with the volume form $dv_g|_{R_C}$ and the Hermitian metrics corresponding to $h$ (see (3.1)) on powers of $ \Cal L _R$. We, then, consider
$$
\Lambda_R:=\lim _{t \to \infty}{1\over t}g_{R,t}. 
$$

We claim it will be a real analytic $Gal(R/C)$-invariant metric on $R_C$. Let  $D_{R,t}$ denote the functional element of  diastasis of $g_{R,t}$ at an arbitrary point of $R_C$.  Since   $L_R^t \subseteq \eusm K^{mt}$ for $\forall t$,               $D_{R,t}$ is bounded by (a multiple of)  $D_{b,R,m,t}$.

%The latter follows from the comparison of inner products in (3.3.2). Indeed, for %$m'\gg m$, let $\Cal L_R^{m'}$ be the inverse image of $\Cal L_C^{m'}$  on %$R_C$. The bundle $ \eusm K_R^m \subset \Cal L_R^{m'}$ has two inner %products, namely,  one with $dv_g|_{R_C}$ and the Hermitian metric coming %from $\Cal L_R^{m'}$, and the standard one on $\eusm K_R^m$ as in the %classical example (3.3.3). We obtain two (equivalent)  metrics on $R_C$. Both %metrics correspond to the bundle $ \eusm K_R^m$.

It follows the uniform convergence of the corresponding 
%diastasic potentials  viewed as complex 
holomorphic functions on $R_C\times \bar R_C$.  Furthermore, $D_{R,t}$ generates a global function on $R_C$. It follows the diastasic potential of the metric $\Lambda_R$ on $R_C$, $\bold P_R:=\bold P_{R,\bold a}$, is a {\it global function} on $R_C$, where $\bold a \in R_C $ is the image of the origin $0\in \Delta$.

Therefore we have established

\proclaim{Proposition-Definition 2} Let
$\Lambda_R :=\lim_{t \to \infty}{1\over t} g_{R,t}$
be a metric  on $R_C$.
It is a real analytic $Gal(R/C)$-invariant Kahler metric. Its diastasic potential
is a function $\bold P_R:=\bold P_{R,\bold a}$ on $R_C$, where $\bold a\in R_C$ is the image of the origin $0\in \Delta$. 
\endproclaim

$(4.2.3)$ {\it  Remark.}\/
 Now, we assume $\dim X \geq 2$, $\eusm K_X$ is ample, $\pi_1(X)$ is residually finite and large, and the genus of the general curvilinear section $C\subset X$ is at least 2. Then $U_X$ may not have the Bergman metric,  as in an example suggested by Campana, namely, consider a sufficiently general ample divisor $X$ in a simple Abelian variety of 
dimension at least  3.

\smallskip
(4.3) Now, we return to the situation in (3.1) with $U=U_X$. For each positive integer $m_0$, the Hermitian metric $h$ on $\Cal L_X$ induces a Hermitian metric $h^{m_0}$ on $\Cal L_X^{m_0}$ as well as on all inverse images of  
$\Cal L_X^{m_0}$ on the coverings of $X$. 

\smallskip
(4.3.1)
	We choose an orthonormal basis $(s^{m_0}_0, \dots, s^{m_0}_{r_{m_0}})$ of   $H^0(X, \Cal L^{m_0}_X)$ with respect to $dv_g$ and $h^{m_0}$.  
We have an inner product and a natural embedding:
$$
\langle s^{m_0}_\alpha, s^{m_0}_\beta \rangle := %{1\over {Vol_g(X)}}
\int_X h^{m_0}(s^{m_0}_\alpha, s^{m_0}_\beta)dv_g; \quad  \phi_{X,{m_0}}: X\hookrightarrow \bold  P^{r_{m_0}}:=\bold  P(H^0(X, \Cal L^{m_0}_X)^*).
$$
Let $g_{FS}$ denote the corresponding standard Fubini-Study metric on the projective space. As in Yau \cite{Y, Sect.\; 6, p.\; 139} (see also Tian  \cite{Ti}), the ${1\over {m_0}}$-multiple of  $g_{FS}$ on $\bold P^{r_{m_0}}$ restricts to a Kahler metric  on  $X$:
$$
g_{X,{m_0}}:={1\over {m_0}} \phi^*_{X,{m_0}}g_{FS}.
$$
\smallskip

(4.3.2) Similar statement holds for all finite coverings of $X$. The bundles $\tau_{i0}^*\Cal L_X$  $(0\leq i<\infty)$ are ample. However, $\tau_{i0}^*\Cal L_X$'s are not necessary very ample  bundles.  

For an appropriate $m_i$, the bundle $(\tau_{i0}^*\Cal L_X)^{m_0m_i}$ is very ample hence it defines a natural embedding
$
\phi_{X_{i},{m_0m_i}}: X_i \hookrightarrow \bold P^{r_{m_0m_i}}
$ into an appropriate projective space.
As above, we get a metric  $g_{X_{i},{m_0m_i}}:= {1\over {m_0m_i}} \phi^*_{X,{m_0m_i}}g_{FS}$  on $X_i$ and the corresponding diastasic potential.  

For all $i=0, \dots$, we consider the metrics $g_{X_{i},{m_0m_i}}$  as $m_0 \rightarrow \infty$.
%Finally, we consider the integers $m_{ij}:=m_i + j$ for $0\leq i, j < \infty$. As %above, we obtain the metrics $g_{X_{i}, m_{ij}}$ on $X_i$.
\smallskip

(4.3.3)  We consider pullbacks on $U_X$  of the  metrics $g_{X_{i},{m_0m_i}}$  and the corresponding diastasises. We will establish that the functional elements of the diastasises converge at a point $\bold p \in U_X$, and we will obtain a real analytic strictly plurisubharmonic functional element at $\bold p$. These functional elements will define the desired Kahler metric $\Lambda$ on $U_X$. 

Let $H_\infty$ denote the hyperplane at infinity in $\bold P^r$.  We can and will assume that  the functional elements of diastasises generate functions on the preimages of $X\backslash H_\infty$ on $X_i$'s.

	We assume the point $\bold p$ does not lie at infinity.  We consider a small compact neighborhood $G\subset U_X$ of $\bold p$. The pullbacks  on $U_X$ of the diastasises are functions on $G$. 
%(As in $\bold P^r \backslash H_\infty$ \cite{C, p.\;17}, it is often convenient to %replace the metric given by $D(p, q)$ in the complement to infinity by the metric %given by $D'(p, q) := e^{D(p, q)}-1$.)

First, we establish the pointwise convergence at $\bold p$ of the pullbacks of  diastasises. We take a sufficiently general curvilinear section $C\subset X\subset \bold P^r$ whose inverse image on $U_X$ contains $\bold p$.  We, then, apply Proposition-Definition 2. 

 Next, we apply the Montel theorem to get the uniform convergence on $G\times \bar G$ of the corresponding holomorphic functional elements  (the complexifications as in Sect.\;2) to a holomorphic functional element. We obtain a {\it real}\/ analytic functional element at the fixed point $\bold p$, {\it denoted}\/ by $D_{U_X}(\bold p,u)$.

 It is easy to see that $D_{U_X}(\bold p, z(u), \bar z(u))$ is strictly plurisubharmonic, where $z, \bar z$ are coordinates in a neighborhood with origin at $\bold p$. Indeed, we take an arbitrary tangent vector $\bold v$ to $U_X$ at $\bold p$. We can assume that the inverse image of a sufficiently general curvilinear section of  $X$ is tangent to $\bold v$. 
We get $D_{U_X}(\bold p,u)$ is strictly plurisubharmonic. We observe that $D_{U_X}(\bold p,u)$ will be the diastasic potential at $\bold p$ of the desired metric $\Lambda$.
% Finally, by the fundamental property of the diastasis, the metric $\Lambda$ is %independent of the very ample bundle $\Cal L_X$.

Thus, we have established the following

\proclaim{Proposition-Definition 3} We assume that a general curvilinear section $C\subset X$ has $g(C)\geq 2$. Then $U_X$ is equipped with a real analytic $\pi_1(X)$-invariant Kahler metric, denoted by $\Lambda$.
 The restriction of $\Lambda$ on $R_C$, the inverse image on $U_X$ of a general curvilinear section $C\subset X$, is the metric $\Lambda_R$ on $R_C$. 
%Such a metric $\Lambda$ is unique, i.e., it is independent of $\Cal L_X$.
\endproclaim

\head
5. Proof of the theorem
\endhead

We assume  $U_X$ and all $X_i$'s are equipped with the metric $\Lambda$. Let $I, R\subset U_X$ be two subsets with $I$ compact. Let  $d(u,p)$ denote the distance function on  $U_X$ with its Riemannian structure induced by $\Lambda$.
 We set
$
  \xi(I,R):=\sup_{u\in I} [\inf_{p\in R} d(u,p)].
$

 We say that a sequence of subsets $\{R_\gamma\}_{\gamma\in \bold N}$, where $R_\gamma\subset U_X$, approximates the set $I$ if
$
\lim_{\gamma \to \infty}\xi(I,R_\gamma) =0.
$
\smallskip
%If $\{R_j\}_{j\in \bold N}$  are, in addition, connected holomorphic subsets that %approximate $I$ then $\lim_{j \to \infty}R_j$ is not necessary a holomorphic subset. In %fact, $I$ then $\lim_{j \to \infty}R_j$ may not exist at all.  Here, we will employ that %$\pi_1(X)$ is large. 

(5.1)  {\it Prolongation}.\/ 
Let $z=(z_1, \dots, z_n)$ be a local coordinate system in a small neighborhood $\Cal V$ with origin at a {\it  fixed} point $\bold a\in \Cal V\subset U_X$. Let $\tilde \Phi_\bold a (z(p), \overline{z(p)})$ be the diastasic potential at $\bold a$. 
Let  $\bold b\in U_X$ be an arbitrary point. Let  
$$
\qquad \qquad\qquad \qquad I: u=u(s) \qquad (0\leq s\leq 1,\; u(0)= \bold a, u(1)= \bold b)
$$
 be a path joining $\bold a$ and $ \bold b$. 
\smallskip

(5.1.1)  {\it Prolongation along the path} $I$.\/ We say $\tilde \Phi_\bold a (z(p), \overline{z(p)})$ has a prolongation along $I$  if the following two conditions are satisfied:

\roster
\item"{(i)}" To every $s\in [0,1]$ there corresponds a functional element of diastasic potential $\tilde \Phi_{u(s)} (z(p), \overline{z(p)})$ at $u(s)$ ($u(s)$ is also called the center).

\item"{(ii)}" For every $s_0\in [0,1]$, we can take a suitable subarc $u:=u(s) \; (|s-s_0|\leq \epsilon, \epsilon >0)$ of $I$ contained in the domain of convergence of  $\tilde \Phi_{u(s_0)} (z(p), \overline{z(p)})$ such that every functional element $\tilde \Phi_{u(s)} (z(p), \overline{z(p)})$ with $|s-s_0|\leq \epsilon$ is a direct prolongation of $\tilde \Phi_{u(s_0)} (z(p), \overline{z(p)})$.
\endroster

The {\it direct prolongation} means the following. Suppose $\tilde \Phi_{\alpha_1} (z(p), \overline{z(p)})$ is defined on $\Cal V_1$ and $\tilde \Phi_{\alpha_2} (z(p), \overline{z(p)})$ is defined on $\Cal V_2$ ($\Cal V_1 \cap \Cal V_2 \not= \emptyset$). Then $\tilde \Phi_{\alpha_2} (z(p), \overline{z(p)})$ is the direct prolongation of $\tilde \Phi_{\alpha_1} (z(p), \overline{z(p)})$ if they coincide on $\Cal V_1 \cap \Cal V_2$.
Recall that the complexification allows us to consider the corresponding holomorphic function in place of the diastasic potential. It follows  a prolongation along $I$ is unique provided it exists. 

Let  $A = \{a_\nu\}$ ($\bold a, \bold b \in A$, $I=\bar A$) be a countable ordered dense subset of $I$. 
We would like to prolongate $\tilde \Phi_\bold a (z(p), \overline{z(p)})$ along $I$ obtaining   the diastasic potential $\tilde \Phi_{a_\nu} (z(p), \overline{z(p)})$ of $\Lambda$ for each $a_\nu$.
We claim the prolongation along $I$ is possible. 
\smallskip

(5.1.2) Now, we will make use of the assumption $\pi_1(X)$ is large. We can assume $I$ is embedded in $X$ via $\tau$; otherwise, we could have replaced $X$ by $X_i$ for $i\gg0$. The set $A$ is a union of an increasing sequence of finite ordered subsets: 
$$
A_1 \subset A_2 \subset \cdots \subset A_\gamma \subset \cdots \subset A, \qquad \bold a, \bold b \in A_\gamma\; (\forall \gamma).
$$
%By the definition of prolongation along a path, the set of points $q\in I$ such that $\tilde %\Phi_\bold a (z(p), \overline{z(p)})$ can be prolongated to $q$ is an open subset of $I$. %We will show this subset is closed.
We consider an arbitrary   $A_\gamma $ and the corresponding  set 
$$
\tau_i(A_\gamma)\subset \tau_i(I)\subset X_i \subset \bold P^{r_i},
$$ 
where $i$ is a sufficiently large integer and $r_i$ is an appropriate integer.  We apply Bertini's theorems to the linear system of curvilinear sections passing through $\tau_i(A_\gamma)$, i.e., the moving part of the system is a one-dimensional subscheme in $X_i$. 
We claim this linear system (and its inverse images) have no fixed components on $X_i$ for all $i\gg 0$.

Suppose, to the contrary, $W\subset X_i$ is a fixed  component. Then $W$ belongs to the linear span of $\tau_i(A_\gamma)\subset \bold P^{r_i}$. We move up along the tower (3.1.1).  For $j\gg i \gg 0$, the linear span of $\tau_j(A_\gamma)\subset \bold P^{r_j}$ will not contain $\tau^{-1}_{ji}(W)$. Hence the corresponding linear system on $X_j$ does not contain $\tau^{-1}_{ji}(W)$.  Therefore the linear system has no fixed components.

A priory, a general member of the system may have singularities at the base points  of  the system. However, we can always assume $\bold P^{r_i}$ is sufficiently large, and our system is sufficiently large as well.
Thus, the general member of the linear system on $X_j$ will be a connected nonsingular curve. Its inverse image on $U_X$ will be a connected open Riemann surface $R_\gamma$ by the Campana-Deligne theorem \cite{Kol, Theorem 2.14}.   These Riemann surfaces will approximate  $I$ as $\gamma$  goes to infinity.
\smallskip

(5.1.3)  
For every $\gamma$, the diastasic potential  of the induced metric on $R:= R_\gamma$ is the restriction of the corresponding diastasic potential of $U_X$, and  $\bold P_R:=\bold P_{R,\bold a}$ is a function on $R$.  

We replace the path $I$ and an arbitrary $A_\gamma\; (\gamma \gg 0)$ by a broken geodesic $\sigma_\gamma$ between the points $\bold a $ and $\bold b$. Namely, we replace the subpath of $I$ between two adjacent points of $A_\gamma$ by a geodesic on $U_X$. 
We also consider the corresponding broken geodesic $\rho_\gamma$ on $R_\gamma$.  Recall  (Sect.\; 2) that the diastasis approximates the square of the geodesic distance in the {\it small}. Hence $\rho_\gamma$ will be close to $\sigma_\gamma$ provided each $a_\nu$ is close to $a_{\nu+1}$,
and we get
$$
\lim_{\gamma \to \infty} \xi(I, \sigma_\gamma) =  \lim_{\gamma \to \infty} \xi(I, \rho_\gamma) = 0 \qquad (\sigma_\gamma \subset U_X, \; \rho_\gamma \subset R_\gamma).
$$
\smallskip

(5.1.4) Now, we will establish the prolongation along $I$. Assume we can prolongate along $I\backslash \bold b$. We take a sufficiently small subarc $E\subset I$ in 
the domain  $\Cal V_\bold b$ of $\tilde \Phi^{\Cal V_\bold b}_\bold b (z(p), \overline{z(p)}) := \tilde \Phi_\bold b (z(p), \overline{z(p)}) $. Take a point   $\bold w\in E\backslash \bold b$ and its small neighborhood $\Cal V_\bold w \subset \Cal V_\bold b$ in $U_X$. We set  $\tilde \Phi^{\Cal V_\bold w}_\bold b (z(p), \overline{z(p)}):= \tilde \Phi^{\Cal V_\bold b}_\bold b (z(p), \overline{z(p)})|\Cal V_\bold w$, more precisely,  $ \tilde \Phi^{\Cal V_\bold w}_\bold b (z(p), \overline{z(p)}) $ is a  real analytic function on $\Cal V_\bold w$ with center $\bold w$ (functional element)  obtained from the real analytic function $ \tilde \Phi^{\Cal V_\bold b}_\bold b (z(p), \overline{z(p)})$ on $\Cal V_\bold b$ with center $\bold b$.
For $p\in \Cal V_\bold w$, we claim
$$
\Phi_{\bold w\bold b}(p) := \tilde \Phi^{\Cal V_\bold w}_\bold b (z(p), \overline{z(p)}) - \tilde \Phi^{\Cal V_\bold w}_\bold w (z(p), \overline{z(p)})= D_{U_X}^{\Cal V_\bold w}(\bold b, p) - D_{U_X}(\bold w,p) =0, 
$$
where $D_{U_X}^{\Cal V_\bold w}(\bold b, p)$ is the real analytic function in $z(p), \overline{z(p)}$ on $\Cal V_\bold w$ with center $\bold w$ obtained from the real analytic function $D_{U_X}(\bold b, p)$ in $z(p), \overline{z(p)}$ on $\Cal V_\bold b$ with center $\bold b$.

Let $u\in \Cal V_\bold w$ be an arbitrary point.  We choose $\{R_\gamma\}$, as above, with an additional condition: $u, \bold w \in R_\gamma  (\forall \gamma)$. Then
$\Phi_{\bold w\bold b}\big|\!R_\gamma(u) = 0$ for all $\gamma \gg 0$ because
$$
D^{\Cal V_\bold w}_{R_\gamma}(\bold b,u) - D^{\Cal V_\bold w}_{R_\gamma}(\bold w,u)=0 \qquad (\forall \gamma \gg 0)
$$
and the fundamental property of diastasis (Sect.\; 2).

It follows we can prolongate $\tilde \Phi_\bold a (z(p), \bar z(p))$ along  $I$.
Since $U_X$ is simply connected, we obtain the desired function $\bold P_U:=\bold P_{U,\bold a}$ on $U_X$.

\smallskip

(5.2) In view of the Oka-Grauert-Narasimhan theorem (Grauert's version), it remains to verify that, for any real $\alpha$, the following set is relatively compact in $U_X$:
$$
E_\alpha :=\{u \in U_X \big |\; \bold P_U(u) < \alpha\}.
$$

Suppose $S\subset E_\alpha$ is an {\it infinite}\/ discrete subset without limit points in $U_X.$ Then we will derive a contradiction by showing that $\bold P_U$ is unbounded on $S$. Since $X$ is compact, $\tau(S)$ will be either a finite set or it will have a limit point. It suffices to replace $S$ by an {\it infinite}\/ set  $T_\alpha$ in the fiber of $\tau$ over a point $Q\in X$ and  
 show that $\bold P_U$ is unbounded on $T_\alpha$. If $\tau(S)$ has a limit point then $Q$ is such a point.

We consider a general curvilinear section $C\subset X$ through $Q$. Set
$
R_C:=\tau^{-1}(C) .
$ 
We obtain a connected open Riemann surface by the Campana-Deligne theorem \cite{Kol, Theorem 2.14}. 
          
By the fundamental property of diastasis (Sect.\;2), $\bold P_R = \bold P _U|_{R_C}$ where $\bold P_R$ is the corresponding diastasic  potentials on $R_C$. 
 One can find an infinite discrete subset $\tilde T_\alpha \subset \Cal F$ whose image on $R_C$ will be close to the corresponding points of  $T_\alpha$. Moreover, $\tilde T_\alpha$ is approaching the boundary of $\Delta.$

We pick a point $\bold q \in \Cal F$ such that its image on $X$ is close to $Q.$ We will identify $D_{\Lambda_R}(\bold q, \cdot)$ with its inverse image on $\Cal F$.
We see that $D_{\Lambda_R}(\bold q, \cdot)$ goes to infinity on $\tilde T_\alpha$
by considering the tower of coverings: 
$$
C\leftarrow \cdots \leftarrow C_i \leftarrow \cdots \leftarrow R_C,  
$$
 where $C_i \subset X_i$ (see (3.1.1)), and the diatasises of the corresponding  Bergman-type metrics of members of the tower restricted to the complement of hyperplane at infinity where they generate functions  as in (4.3.3). 

By Proposition 1, the diastasises increase as we move up in the tower. 
  So, $\bold P_R$ is unbounded on $\tilde T_\alpha$ and $T_\alpha,$ and $\bold P_U$ will be unbounded on $S$.

The contradiction proves the theorem.
\medskip
$(5.3)$ {\it  Remark.}\/ Similar argument proves the Shafarevich conjecture when $X$ is singular. Namely, let $X \hookrightarrow  \bold P^{r}$ be a  connected normal projective variety of dimension $n>0$. Assume $\pi_1(X)$ is large and residually finite. Then its universal covering is a Stein space. The details will appear elsewhere.

%Presumably, in the above, one can replace {\it projective variety}\/ by a normal %{\it Kahler}\/ space $X$.

\medskip
{\it  Acknowledgment.} The author would like to thank Joseph Bernstein, J\'anos Koll\'ar, and Pierre Milman for helpful conversations and the referee for helpful remarks. The author is grateful to Fedor Bogomolov, Fr\'ed\'eric Campana,  Takeo Ohsawa, and Peter Polyakov for their  emails.

\enddocument

%\ref  \key GH  \by P. Griffiths, J. Harris
%\book Principles of Algebraic Geometry
%\publ John Wily, New York
%\yr 1978
%\endref

%\ref \key Kr \by 	I. Kra \book  Automorphic forms and Kleinian groups
%\publ Mathematics Lecture Note Series, W. A. Benjamin, Inc., Reading, Mass. 
%\yr 1972 \endref

%\ref  \key L  \by  J. Lehner \book Discontinuous Groups and Automorphic Functions
%\publ American Math. Soc., Providence \yr 1964 \endref

%\ref  \key Ma  \by  A. Mal'cev \pages  405--422
%\paper On isomorphic matrix representations of infinite groups
%\yr1940 
%\vol  8(50)
%\jour Mat. Sbornik
%\endref

%\ref  \key D  \by  D. Drasin \pages 356--365
%\paper Cusp forms and Poincar\'e series
%\yr1968 \vol  90
%\jour Amer. J.  Math.
%\endref

%\ref 
%\key E  \by  P.  Eyssidieux 
%\paper  Lectures on Shafarevich conjecture on uniformization 
%\yr 2011
%\publ Inst. Fourier,  Grenoble (France).
%($ www.fourier.ujf\;grenoble.fr/eyssi/cours\;eyssidieux \;conjecture$ %$\;de\;shafarevich.\;v4.pdf
%$)
%\endref

%\ref  \key  EKPR  \by   P. Eyssidieux, L. Katzarkov, T. Pantev, M. Ramachandran
% \pages  1545 --1581
%\paper Linear Shafarevich conjecture
%\yr 2012 \vol  176
%\jour Ann. of Math.
%\endref

JANUARY 17, 2014

(4.2.2) Now, let $\Cal L _R$ denote the inverse image of $\Cal L_C$ on $R_C$. We consider the metrics $g_{R,t}$ constructed, as in \cite{Ti}, with the volume form $dv_g|_{R_C}$ and the Hermitian metrics corresponding to $h$ (see (3.1)) on powers of $ \Cal L _R$. We, then, consider
$$
\Lambda_R:=\lim _{t \to \infty}{1\over t}g_{R,t}. 
$$

We claim it will be a real analytic $Gal(R/C)$-invariant metric on $R_C$. Let  $D_{R,t}$ denote the functional element of  diastasis of $g_{R,t}$ at an arbitrary point of 
$R_C$.  It is bounded by (a multiple of)  $D_{b,R,m,t}$ considered in (4.2.1). The latter follows from the comparison of inner products in (3.3.2). Indeed, for $m'\gg m$, let $\Cal L_R^{m'}$ be the inverse image of $\Cal L_C^{m'}$  on $R_C$. The bundle $ \eusm K_R^m \subset \Cal L_R^{m'}$ has two inner products, namely,  one with $dv_g|_{R_C}$ and the Hermitian metric coming from $\Cal L_R^{m'}$, and the standard one on $\eusm K_R^m$ as in the classical example (3.3.3). We obtain two (equivalent)  metrics on $R_C$. Both metrics correspond to the bundle $ \eusm K_R^m$.

JANUARY 2, 2014

(3.3.5) {\it Comparison of inner products.}\/ Given the volume form $dv_g$ and the Hermitian metric $h$ on $\Cal L = \eusm K^t_X$, one can define an inner product on sections of $\Cal L$ in two ways:
$$
\langle \omega_1,\bar \omega_2\rangle' := \int_U h(\omega_1,\bar\omega_2)dv_g
\qquad  {\text {and}} \qquad
\langle \omega_1,\bar \omega_2\rangle'' := \int_U \Cal H_{U,t}(\omega_1\otimes \bar\omega_2).
$$
As in \cite{Kol, Chap.\;5.1}, one can compare the corresponding two norms. 

We take a hyperplane $H_\infty \subset \bold P^r$ and consider a point $\bold u\in U$ in the inverse image of $X\backslash H_\infty$ on $U$. We consider local coordinates $z_i$ at $\bold u$.
Let
$$
\omega:= g(\bold u) (dz_1\wedge \cdots \wedge dz_n)^t,
$$
$$
h_{\eusm K^t_U}(\bold u): =h^t_{\eusm K_U}\bold (\bold u):=  \|dz_1\wedge\cdots\wedge dz_n \|^t, \;   {\text {where}} \; \|\cdot\|= h^{1/2}_{\eusm K_U}(.,.),   \quad                              {\text {and}}
$$
$$
dv:=dv_g=V_g\! \prod^n_{i=1}(\sqrt{-1}\; dz_i\!\wedge\! d\bar {z}_i).
$$
Then
$$
h_{\eusm K_U^t}(\omega)^2dv_g =  h^{2t}_{\eusm K_U}(\bold u)|g(\bold u)|^2\!\cdot\! V_g\! \prod^n_{i=1}(\sqrt{-1}\; dz_i\!\wedge\! d\bar {z}_i), \quad                                       {\text {and}}
$$
$$
\Cal H_{U,t}(\omega,\bar \omega)= |g(\bold u)|^{2/t} dz_1\wedge \cdots \wedge dz_n \wedge d\bar z_1\wedge \cdots \wedge d \bar z_n.
$$
\smallskip

November 27--2013

(3.2.3) The function $B(z,w)$ can be replaced  by a section of a relevant bundle as in  \cite{Kol, Chap.\;7, pp.\;81-84, Lemma-Definition 7.2}.
We consider
$M\times M$ and write a point of 
$M\times M$ as a pair
$(z, w)$.  The reproducing kernel, provided it exists, will be a section of a bundle $ p^*_1 \Cal L \otimes p^*_2\bar \Cal L$  where, now, $\Cal L$ is a line bundle on $M$ and $p_1$ and $p_2$ are the coordinate projections of $ M\times M.$ 
Let
$\{u_k\}\subset H\subset H^0(M,\Cal L)$ denote an orthonormal basis of some Hilbert space $H.$ Set
$$
B(z,w):=\sum_k p^*_1 u_k(z)\otimes p^*_2 \bar u_k(w).
$$
If we assume the evaluation map at every point $Q\in M$ is a continuous linear functional
then we get a map to a projective space as above (see also \cite{Kob2, Chap.\;4.10, pp.\;224-228}). In (4.3.3) below, we consider the case when, in general, we do not get the natural map to a projective space.  However, we will be able to define  $B(z,w)$ locally. 

(3.2.4) {\it Projective spaces and Veronese maps}.\/ For $p_0 \in\bold P^r$, we consider canonical coordinates $(z_1, \dots, z_r)$ with center $p_0$ on the complement of a hyperplane at infinity.  By Calabi (\cite{Chap.\;4, (27)}, Calabi considered $\bold P^r$ with $r=\infty$ as well),
$$
	D_{\bold P^r}(p_0,p)= \log\bigl(1+\sum_\sigma|z_\sigma(p)|^2 \bigl).
$$
In the homogeneous coordinates $\xi_0,\dots,\xi_r  $, where $z_\sigma := {\xi_\sigma/\xi_0}$, we get
$$
D_{\bold P^r}(p_0,p)= \log {\sum^r_{\sigma=0} |\xi_\sigma(p)|^2 \over |\xi_0(p)|^2}.
$$

 The tautological bundle $\Cal O_{\bold P^r} (1)$ defines the tautological polarization. The bundle $\Cal O_{\bold P^r} (t)$ defines the Veronese embedding. We get the corresponding Hermitian metrics $h_{\Cal O_{\bold P^r} (1)}$ and $ h_{\Cal O_{\bold P^r} (t)}=h^t_{\Cal O_{\bold P^r} (1)} $, and polarizations on $\bold P^r$. The diastasic potential is not a function on $\bold P^r$ while the diastasic potential of the Poincar\'e metric on  $\Delta$ is a function on $\Delta$ because the bundle $\Cal O_{\bold P^r}(1)$ is not trivial while $\eusm K_\Delta$ is trivial.
\smallskip

May 9

For a given $A_\gamma$, we consider $\tau(A_\gamma) \subset X_i$  $(i\gg 0)$. We consider a sufficiently  large linear system $L_i$ of curvilinear sections on $X_i \subset \bold P^{r_i}$ through $\tau(A_\gamma)$, i.e., the movable part of $L_i$ consists of curves. 

If $L_i$ has {\it no}\/ fixed components then its general member is a nonsingular connected curve $C\subset X_i$
 by  Bertini's and Lefschetz's theorems (compare [GH, pp.\;173-174]). Further, each component of $\tau_i^{-1}(C) \subset U_X$ is an open Riemann surface because $\pi_1(X)$ is large, and the points of $A_\gamma$ lie on the same (connected) Riemann surface, say $R_\gamma$. 
By the Campana-Deligne theorem \cite {Kol, Theorem 2.14} and our assumption, $g(C) \geq 2$. 

We consider the restriction of the metric $\Lambda$ on $R_C$. Consider a broken geodesic on $R_C$ that will replace $I$ and the points $A_\gamma$, where the new real curve $[a_\nu, a_{\nu+1}]\subset R_C$ is a geodesic. Recall that $a_\nu$ was close to $a_{\nu+1} $ 
and  the diastasis approximates the square of the geodesic distance in the {\it small}
\cite {C, Chap.\;2, Prop.\;5}.

The diastasic potential of the open Riemann surface $R_C$, with induced metric, is a function. In fact, it is equal to $\lim_{t \to \infty}{1\over t} \log \bold B_{R,\eusm K^t}(z.\bar z)$, where $\eusm K_{R_C}$ is trivial  because any holomorphic bundle is trivial on an open Riemann surface.  Hence we can prolongate $\tilde \Phi_\bold a (z(p), \overline{z(p)})$ along the broken geodesic on $R_C$.

Now, we will show that $L_i \subset X_i$ has no fixed components provided $i\gg 0$. Given a finite set $A_\gamma$, the possible fixed components will disappear when we move up in the tower (3.1.1) because $\pi_1(X)$ is large.

________
(5.1)
We assume our $U_X$ and all $X_i$'s are equipped with the metric $\Lambda$. Let $I\subset U$ be an arbitrary compact subset and $R\subset U_X$ an arbitrary  subset. In the sequel,  $I$ will be a path and $R$ will be an analytic subset. For $u\in I$ and $p\in R$, we set
$$
d(u,R):=\inf_{p\in R} d(u,p), \qquad  \xi(I,R):=\sup_{u\in I} d(u,R)
$$
where  $d(u,p)$ is the distance function on the Riemann manifold $U_X$. We say that a sequence of subsets $\{R_j\}_{j\in \bold N}$ approximates  $I$ if
$$
\lim_{j \to \infty}\xi(I,R_j) =0.
$$
Let $\bold a$ and $\bold b$ are two points on $U_X$ and let  
$$
I: u=u(s) \; (0\leq s\leq 1, u(0)= \bold a, u(1)= \bold b)
$$
 be a path joining $\bold a$ and $ \bold b$. 
%By abuse of notation, we denote by $I$ its image in $U_X$.
Let  $A = \{a_\nu\}$ ($I=\bar A$) be an a countable dense set on the path. 
We would like to prolongate $\tilde \Phi_\bold a (z(p), \overline{z(p)})$ along $I$ obtaining   the diastasic potential $\tilde \Phi_{a_\nu} (z(p), \overline{z(p)})$ of $\Lambda$ for each $a_\nu$.
We claim the prolongation along $I$ is possible.

The set $A$ is a union of increasing sequence of finite subsets: 
$$
A_1 \subset A_2 \subset \cdots \subset A_\gamma \subset \cdots \subset A, \qquad \bold a, \bold b \in A_\gamma\; (\forall \gamma).
$$
 For a given $A_\gamma$, we consider $\tau(A_\gamma) \subset X_i$  $(i\gg 0)$. We consider a sufficiently  large linear system $L_i$ of curvilinear sections on $X_i \subset \bold P^{r_i}$ through $\tau(A_\gamma)$, i.e., the movable part of $L_i$ consists of curves.

May 31

The bundles $\Cal L^t$ ($t$ is a positive integer) define the embeddings of $X$ in projective spaces. We consider $\tau(A_i)$. For a sufficiently large $t_i$, we consider a general curvilinear section of $X$ through $\tau(A_i)$ in the corresponding projective space. We obtain a sequence $\{C_i\}$ of nonsingular curves on $X$ where $A_i \subset C_i$.

We claim $\Cal C:=\lim_{i \to \infty} C_i$ is a holomorphic curve, and $\tau ^{-1}(\Cal C_{\text{red}}) \subset U_X$ is a connected open one-dimensional Stein subsubspace containing our path $I$. Indeed, $\Cal C$ is a closed subset of $X$ and, in a small neighborhood of each point of $X$, it is a zero-set of a holomorphic function. Moreover, $\tau ^{-1}(\Cal C_{\text{red}})$ is open because $\pi_1(X)$ is large.

May 27

We consider
a path $I =[\bold a, \bold b]\subset U_X$ and a finite number of points $\bold a =a_1, a_2, \dots, a_\kappa=\bold b$ ($a_\nu \in I$, $\kappa < \infty$) where $a_\nu$ is close to $a_{\nu+1}$.  
We would like to prolongate $\tilde \Phi_\bold a (z(p), \overline{z(p)})$ along $I$ obtaining   the diastasic potential $\tilde \Phi_{a_\nu} (z(p), \overline{z(p)})$ of $\Lambda$ for each $a_\nu$.
We claim the prolongation along $I$ is possible.

We consider $\{\tau_i(a_\nu)\} \subset X_i (i\gg 0)$. We consider a pencil $L$ of curves on $X_i$ with
base points $\{\tau_i(a_\nu)\}$. We can assume its general member is a nonsingular connected curve $C\subset X_i$
 by  Bertini's and Lefschetz's theorems (compare [GH, pp.\;173-174]). Further, each component of $\tau_i^{-1}(C) \subset U_X$ is an open Riemann surface because $\pi_1(X)$ is large. By the Campana-Deligne theorem \cite {Kol, Theorem 2.14} and our assumption, $g(C) \geq 2$. 

We claim that the point $\{a_\nu\}$ belong to the same component, say $R_C$. We consider an open exhaustion of $U_X$: $\{U_r \subset U_X, r=1,2, \dot\}$ where each $\bar U_r$ is compact and $\bar U_r \subset U_{r+1}$. If $a'$ and $a''$ belong to two different component then their images in $X_j$, for $j\geq 0$, belong to two different components of $\tau^{-1}_{j0}(C)$. However,  $\tau^{-1}_{j0}(C)$ is connected
by the Campana-Deligne theorem \cite {Kol, Theorem 2.14}.

We consider the restriction of the metric $\Lambda$ on $R_C$. We consider a broken geodesic in $R_C$ that will replace our path $I$, where the new real curve $[a_\nu, a_{\nu+1}]\subset R_C$ is a geodesic; recall that $a_\nu$ was close to $a_{\nu+1} $ (see \cite {C, Chap.\;2, Prop.\;5} or Sect. 2).

The diastasic potential of the open Riemann surface $R_C$, with induced metric, is a function. In fact, it is equal to $\lim_{t \to \infty}{1\over t} \log \bold B_{R,\eusm K^t}(z.\bar z)$, where $\eusm K_{R_C}$ is trivial  because any holomorphic bundle is trivial on an open Riemann surface.  Hence we can prolongate $\tilde \Phi_\bold a (z(p), \overline{z(p)})$ along the broken geodesic on $R_C$.

Since the paths like $I$ can be approximated by broken geodesics
and $U_X$ is simply connected,  we can prolongate along every $I$ and we obtain the desired function $\bold P_U$ on $U_X$.

May 22  Sect. 5

We consider $\{\tau_i(a_\nu)\} \subset X_i (i\gg 0)$. Take a general curvilinear section of $X$ through the points $\{\tau(a_\nu)\}$, i.e., this section contains a 1-dimensional irreducible subvariety as a component (a priory, other components can have dimension $> 1$). If this is impossible then we replace $X$ with $X_i$, $i\gg 0$ (see (3.1.1)). We can assume a linear system $L$ on $X_i$ of curvilinear sections through $\{\tau_i(a_\nu)\}$ has a sufficiently large dimension.

By the Bertini theorem, if  $L$ has no fixed components then its general member is nonsingular and, even, connected by the Lefschetz theorem (compare [GH, pp.\;173-174]). So far, we have not used that $\pi_1(X)$ is large.

We claim if $i\gg 0$ than $L$ has no fixed components. Suppose $W$ is a fixed component
containing the points $\{\tau_i(a_\beta)\} \subset \{\tau_i(a_\nu)\}$ where $\{a_\beta\}$
 are all the points projecting to $W$. Since $\pi_1(X)$ is large, the corresponding component $W' \subset U_X$ will be a compact, where $\tau_i(W')= W$ and we assume, as  we may, $\{a_\beta\} \subset W'$.  Indeed, $W'$ contains no infinite discrete subsets
because, otherwise, a fiber of $\tau_i :W' \rightarrow W$ will contain an infinite discrete subset. this can not happen for $i\gg 0$.

Thus, we can assume that a general member of a large linear system $L$ of curvilinear sections on $X_i$ $(i\gg0)$ is a nonsingular connected curve $C$ of genus $g(C)\geq 2$. Indeed, by the Campana-Deligne theorem \cite {Kol, Theorem 2.14}, $\tau_i^{-1}(C)$ $(i\gg0)$ will be a connected open Riemann surface $R_C$. 
By the assumption (see Sect.\;1), $\pi_1(X)$ is not Abelian. Hence $\pi_1(C)$ is not Abelian.

aaaaaaaaaaaaa

Take a general curvilinear section of $X$ through the points $\{\tau(a_\nu)\}$, i.e., this section contains a 1-dimensional irreducible subvariety as a component (a priory, other components can have dimension $> 1$). If this is impossible then we replace $X$ with $X_i$, $i\gg 0$ (see (3.1.1)). We can assume a linear system $L$ on $X_i$ of curvilinear sections through $\{\tau_i(a_\nu)\}$ has a sufficiently large dimension.

By the Bertini theorem, if  $L$ has no fixed components then its general member is nonsingular and, even, connected by the Lefschetz theorem (compare [GH, pp.\;173-174]). So far, we have not used that $\pi_1(X)$ is large.

We claim if $i\gg 0$ than $L$ has no fixed components. Suppose $W$ is a fixed component
containing the points $\{\tau_i(a_\beta)\} \subset \{\tau_i(a_\nu)\}$ where $\{a_\beta\}$
 are all the points projecting to $W$. Since $\pi_1(X)$ is large, the corresponding component $W' \subset U_X$ will be a compact, where $\tau_i(W')= W$ and we assume, as  we may, $\{a_\beta\} \subset W'$.  Indeed, $W'$ contains no infinite discrete subsets
because, otherwise, a fiber of $\tau_i :W' \rightarrow W$ will contain an infinite discrete subset. this can not happen for $i\gg 0$

Thus, we can assume that a general member of a large linear system $L$ of curvilinear sections on $X_i$ $(i\gg0)$ is a nonsingular connected curve $C$ of genus $g(C)\geq 2$. Indeed, by the Campana-Deligne theorem \cite {Kol, Theorem 2.14}, $\tau_i^{-1}(C)$ $(i\gg0)$ will be a connected open Riemann surface $R_C$. 
By the assumption (see Sect.\;1), $\pi_1(X)$ is not Abelian. Hence $\pi_1(C)$ is not Abelian.

$$
bbbbbbbbbbbbbbbbb
$$

May 21 Sect. 5

Now, if $\Psi \subset X_i$ is a fixed component of $L$ then we move up in the tower (3.1.1). Since $\tau_i^{-1}(\Psi) \subset U_X$ is not compact, the degree of each component 
$$
W\subset \tau_{ji}^{-1}(\Psi) \subset X_j\subset  \bold P^{r_j} \qquad (j\gg i)
$$ 
will be larger than the dimension $\mu$ of the linear span of all those points in the set $\{\tau_j(a_\nu)\}$ that belong to $W$. Indeed, $\deg(W) \geq \mu +1$ unless $W$ is a rational normal curve of degree  $\mu$ in $\bold P^\mu$ \cite{GH, pp.\;173-174, p.\;179} contradicting the assumption that $\pi_1(X)$ is large. Therefore
we can find a hyperplane in $\bold P^{r_j}$ that does not contain $W$ but contains all the points that are projected to $W$.

Similarly, by moving up in the tower (3.1.1), we can get rid of all possible fixed components containing the points that are images of the  set $\{a_\nu\}$.
========================================

(4.2.2)

If $|\bold B_{\Delta, \eusm K^t}(z,\zeta)| <  |\bold B_{\Delta, \eusm K^t}(\alpha  z,\zeta)|$ for a unique element $\alpha\in \Gamma $ than
$$
{1\over t} \log \bold B_{R,\eusm K^t}(z,\bar z) = {1\over t} \log \bold B_{\Delta, \eusm K^t}(\alpha z,\zeta) +
$$
$$ {1\over t}\log \biggl [Jac_\alpha^t (z)+ \sum_{\gamma\in \Gamma,\gamma \not = \alpha} {\bold B_{\Delta,\eusm K^t}(\gamma z,\zeta) \over \bold B_{\Delta, \eusm K^t}(\alpha z,\zeta)}Jac_\gamma^t(z)\biggl].
$$
Taking the limits we get
$$
\lim _{t \to \infty}{1\over t} \log \bold B_{R,\eusm K^t}(z, \bar z)=\lim _{t \to \infty}{1\over t} \log \bold B_{\Delta,\eusm K^t}(\alpha z,\bar z) > \lim _{t \to \infty}{1\over t} \log \bold B_{\Delta,\eusm K^t}(z,\bar z),
$$
a contradiction. Similar argument works if there is only a finite set of element $\alpha \in\Gamma$ such that  $|\bold B_{\Delta, \eusm K^t}(z,\zeta)| <  |\bold B_{\Delta, \eusm K^t}(\alpha  z,\zeta)|$.

AAAAAAAAAAAAAAAAAAAAAAAAAAAAAAAAAAAAAAAAAAAA

Hence, for $z\in \Cal F$,
$$
\lim_{t \to \infty}{1\over t} g_{R,\eusm K^t}(z, \bar z)=\lim_{t\to\infty} {1\over t} {\partial^2\log\bold B_{R,\eusm K^t}(z,\bar z)\over \partial z \partial \bar z}
% ={\partial^2\log\bold B_{\Delta,\eusm K}(z,\bar z)\over \partial z \partial %\bar z}
={\partial^2\log\bold B_{\Delta,\eusm K}(z,\bar z)\over \partial z \partial \bar z}.
$$

Indeed, the metric $\lim_{t \to \infty} g_{R,\eusm K^t}$ on $R_C$ can be replaced by $ g_{R,\eusm K^t}$ on $R_C$ for an integer $t\gg 0$. Hence it will suffice to establish that $T_{R,\alpha}$ is finite with the {\it new metric}\/
and its diastasic potential (which is close to the potential of $\lim_{t \to \infty} g_{R,\eusm K^t}$).

We proceed as follows. 
We consider $Q\in X\subset \bold P^r$ as the origin of the canonical coordinate system in $\bold P^r\backslash H_\infty$, i.e., $Q=p_0$ and $H_\infty=\{\xi_0=0\}$ in the notation of (3.2.4). We take the inverse images of $H_\infty$ under all $\tau_i$'s. Since $\pi_1(X)$ is residually finite, we get 
$$
T_{R,\alpha}= \bigcup \tau_i^{-1}(T_{X_i,\alpha}),\quad \text{where}\quad T_{X_i,\alpha}:=\{x\in X_i \, \big | \,\tau_{i0}(x)=Q, \, \bold P_{\tau^{-1}_{i0}(C)}(x) <\alpha\}.
$$

As before, we consider the metric on $\tau^{-1}_{i0}(C)$ corresponding  to $\eusm K^t_{\tau^{-1}_{i0}(C)}$.

Next, we consider 
$ 
T_{\Delta,\alpha}=\{u \in \Delta \; |\; \tau_\Delta(u)=Q,\; \bold P_\Delta(u) < \alpha\},  $
 where $\tau_\Delta : \Delta \rightarrow C$ is the projection and we consider
the metric on $\Delta$ corresponding  to $\eusm K^t_\Delta$.
As above, one can describe the set $T_{\Delta,\alpha}$ via the finite coverings of $C$.

 Since the finite quotients of $Gal(R/C)$ are finite quotients of $Gal(\Delta/C)$ as well, we get Card $T_{R,\alpha}$ $\leq$ Card $T_{\Delta,\alpha} < \infty$, the last inequality being classical (Introduction).

% (Card$|T_{\Delta,\alpha}| < \infty$.

%the above inequality $\bold P_R(u) < \alpha$ is equivalent to $\bold B_{R,\eusm K^t}(u,\bar %u)< t\alpha$ for $t\gg 1$.

As before, we consider the metric on $\tau^{-1}_{i0}(C)$ corresponding  to $\eusm K^t_{\tau^{-1}_{i0}(C)}$.

Next, we consider 
$ 
T_{\Delta,\alpha}=\{u \in \Delta \; |\; \tau_\Delta(u)=Q,\; \bold P_\Delta(u) < \alpha\},  $
 where $\tau_\Delta : \Delta \rightarrow C$ is the projection and we consider
the metric on $\Delta$ corresponding  to $\eusm K^t_\Delta$.
As above, one can describe the set $T_{\Delta,\alpha}$ via the finite coverings of $C$.

 Since the finite quotients of $Gal(R/C)$ are finite quotients of $Gal(\Delta/C)$ as well, we get Card $T_{R,\alpha}$ $\leq$ Card $T_{\Delta,\alpha} < \infty$, the last inequality being classical (Introduction).

% (Card$|T_{\Delta,\alpha}| < \infty$.

%the above inequality $\bold P_R(u) < \alpha$ is equivalent to $\bold B_{R,\eusm K^t}(u,\bar %u)< t\alpha$ for $t\gg 1$.

As before, we consider the metric on $\tau^{-1}_{i0}(C)$ corresponding  to $\eusm K^t_{\tau^{-1}_{i0}(C)}$.

Next, we consider 
$ 
T_{\Delta,\alpha}=\{u \in \Delta \; |\; \tau_\Delta(u)=Q,\; \bold P_\Delta(u) < \alpha\},  $
 where $\tau_\Delta : \Delta \rightarrow C$ is the projection and we consider
the metric on $\Delta$ corresponding  to $\eusm K^t_\Delta$.
As above, one can describe the set $T_{\Delta,\alpha}$ via the finite coverings of $C$.

 Since the finite quotients of $Gal(R/C)$ are finite quotients of $Gal(\Delta/C)$ as well, we get Card $T_{R,\alpha}$ $\leq$ Card $T_{\Delta,\alpha} < \infty$, the last inequality being classical (Introduction).

% (Card$|T_{\Delta,\alpha}| < \infty$.

%the above inequality $\bold P_R(u) < \alpha$ is equivalent to $\bold B_{R,\eusm K^t}(u,\bar %u)< t\alpha$ for $t\gg 1$.

As before, we consider the metric on $\tau^{-1}_{i0}(C)$ corresponding  to $\eusm K^t_{\tau^{-1}_{i0}(C)}$.

Next, we consider 
$ 
T_{\Delta,\alpha}=\{u \in \Delta \; |\; \tau_\Delta(u)=Q,\; \bold P_\Delta(u) < \alpha\},  $
 where $\tau_\Delta : \Delta \rightarrow C$ is the projection and we consider
the metric on $\Delta$ corresponding  to $\eusm K^t_\Delta$.
As above, one can describe the set $T_{\Delta,\alpha}$ via the finite coverings of $C$.

 Since the finite quotients of $Gal(R/C)$ are finite quotients of $Gal(\Delta/C)$ as well, we get Card $T_{R,\alpha}$ $\leq$ Card $T_{\Delta,\alpha} < \infty$, the last inequality being classical (Introduction).

% (Card$|T_{\Delta,\alpha}| < \infty$.

%the above inequality $\bold P_R(u) < \alpha$ is equivalent to $\bold B_{R,\eusm K^t}(u,\bar %u)< t\alpha$ for $t\gg 1$.

%% It is easy to see (\cite {C, (23)}) that the diastasis of any Bergman metric equals $\log %%Ker(z,z)$, where $Ker(z,z)$ denotes the corresponding  kernel.

(4.3.1) Now, we return to the situation in (3.1) with $U:=U_X$. For each positive integer $m_0$, the Hermitian metric $h$ on $\Cal L$ induces a Hermitian metric   $h^{m_0}$  on $\Cal L^{m_0}_X$ as well as on all inverse images of $\Cal L^{m_0}_X$ on the coverings of $X$. 

	One chooses an orthonormal basis $(s^{m_0}_0, \dots, s^{m_0}_{r_{m_0}})$ of the space $H^0(X, \Cal L^{m_0}_X)$ of all global sections of $\Cal L^{m_0}_X$. We get an inner product and a natural embedding:
$$
\langle s^{m_0}_\alpha, s^{m_0}_\beta \rangle := {1\over {Vol_g(X)}}\int_X h^{m_0}(s^{m_0}_\alpha, s^{m_0}_\beta)dV_g; \qquad  \phi_{X,{m_0}}: X\hookrightarrow \bold  P^{r_{m_0}}.
$$
Let $g_{FS}$ be the standard Fubini-Study metric on $\bold P^{r_{m_0}}$.  The ${1\over {m_0}}$-multiple of  $g_{FS}$ on $\bold P^{r_{m_0}}$ restricts to a Kahler metric  on  $X$:
$$
g_{X,{m_0}}:={1\over {m_0}} \phi^*_{X,{m_0}}g_{FS}.
$$
One of the main results of Tian in \cite{Ti, Theorems A} (Yau's problem) says 
$g_{X,{m_0}}$ converge, as $m_0 \rightarrow \infty$, to $g_X$ in  $C^2$-topology. 
%His proof is local in nature.

(4.3.2) Similar statement holds for all finite coverings of $X$. The bundles $\tau_i^*\Cal L_X$  $(0\leq i<\infty)$ are ample. However, $\tau_i^*\Cal L_X$'s are not necessary very ample  bundles.  

For an appropriate $m_i$, the bundle $(\tau_i^*\Cal L_X)^{ {m_i}}$ is very ample hence it defines an embedding
$
\phi_{X_{i},{m_i}}: X_i \hookrightarrow \bold P^{r_{m_i}}. 
$
As above, we get a metric  $g_{X_{i},{m_i}}$  on $X_i$ and the corresponding diastasic potential.  

Finally, we consider the integers $m_{ij}:=m_i + j$ for $0\leq i, j < \infty$. As above, we obtain the metrics $g_{X_{i}, m_{ij}}$ on $X_i$.

%As we have mentioned, the volume forms are arising from $X$, and the Hermitian metrics on %the corresponding bundles are arising from the Hermitian metric $h$ on  $\Cal L_X$. The %pullbacks on $\Cal V_p$ of the metrics on $X_i$'s  are {\it smaller}\/ pseudo-metrics than the %corresponding metrics on $\Cal V_p$.

(4.3.4) The corresponding metrics $g_{X_{i}, m_{ij}}$ as well as their pullbacks to $U_X$ have their diastasises. We claim that the corresponding diastasises converge at a point $p \in U_X$ and we obtain a real analytic strictly plurisubharmonic functional element at $p$. Then these functional element will define a Kahler metric on $U_X$, and it will be the desired $\Lambda$.

\proclaim{Proposition-Definition 3} The universal covering $U_X$ is equipped with Kahler metric, denoted by $\Lambda$, such that induced metric on the preimage of every general hyperplane section of $X$ is the Poincar\'e metric. 
\endproclaim

\head

AAAAAAAAAAAAAAAAAAAAAAAAAAAAAAAAAAAAAAAAAAAAAAAAAAAAAAAAAAAAAA

1.1. Given $U, \,\Cal L$ and $B_{U, \Cal L}(z,w)$ as above, we will
show that $U$ has a natural Kahler metric. Since $B(z,z):=B_{U, \Cal L}(z,z) >0,$ we
may consider $\log B(z,z)$.  We will show that $\log B(z,z)$ is {\it strictly}\/ plurisubharmonic,
i.e., 
$$
 \qquad \qquad\qquad\qquad\qquad\qquad  ds^2_U =2 \sum g_{j k} dz_j d\overline z_k
\quad\quad ( g_{j k}:={\partial^2\over\partial z_j \partial \overline z_k  }\log B(z, z) ) 
$$
determines a Hermitian metric.  For an arbitrary real tangent vector  at $z\in U,$  
$$
\qquad\qquad\qquad\qquad\qquad\qquad  \bold v=   \partial_a+ \overline{\partial_a}   
\qquad(\partial_a:=
\sum_{t=1}^n a_t{\partial \over {\partial  z_t}},\quad 
\overline{\partial_a}:=
\sum_{t=1}^n \bar a_t{\partial \over {\partial \bar z_t}}),
$$
we set $|\bold v|^2 :=
\sum g_{jk} a_j{\bar a}_k
$   hence $|\bold v|^2=\partial_a \overline{\partial_a} \log B(z,z) $. The $g_{jk}$'s determine a
Hermitian metric iff $|\bold v|^2>0$ for $\bold v\neq 0$. With the above notation, 
$$ 
\!\!\!\!\!\partial_a \overline{\partial_a} \log B(z,z)= {1\over B(z,z)^2}\Big\{
\sum_k|\bold v\psi_k(z)|^2
\cdot
\sum_k|\psi_k(z)|^2 -\Big|
\sum_k \bold v\psi_k(z)\!\cdot \!\overline {\psi_k(z)}\Big|^2 \Big\}
$$
$$
\qquad  \quad\quad\quad\quad\quad  \geq {1\over B(z,z)^2}\Big\{ \sum_k|\bold v\psi_k(z)|^2
\cdot
\sum_k|\psi_k(z)|^2 -\Big(
\sum_k \big |\bold v\psi_k(z)\!\cdot\! \overline {\psi_k(z)}\big|\Big)^2\Big\}
$$
where $B(z, z)=\sum \psi_k(z) \overline{\psi_k(z)}$ in a neighborhood of $z$ (the calculations
are similar to
\cite {FK, pp.\, 13 --14, pp.\, 190 --191}).
   In the Lagrange  identity
$$
\sum_{k=1}^s A^2_k\,\sum_{k=1}^s B^2_k -\Bigl(\sum_{k=1}^s A_kB_k \Bigl)^2 ={1\over
2}\sum_{i=1}^s \sum_{j=1}^s (A_iB_j-B_iA_j)^2\quad (\forall A_i, \forall B_i \in\bold R_+, \, 
s\in \bold N),
$$
the right-hand side  is strictly positive for $s=\infty,$ provided it is strictly
positive for some $s\in\bold N$ (and the series are convergent). It follows
that $|\bold v|^2 > 0$ for $\bold v\neq 0$, and we get the desired metric.

\ref  
\key D  \by  D. Drasin \pages 356--365
\paper Cusp forms and Poincar\'e series
\yr1968 \vol  90
\jour Amer. J.  Math.
\endref

\ref  
\key  U  
\by M. Umehara \pages  203--214
\paper  Kaehler Submanifolds of Complex Space Forms
\jour Tokyo J. Math
\vol 10, No. 1
\endref
\endRef

\enddocument

(3.4) We pull back all the  metrics in the tower of {\it finite} coverings of $X$ to the neighborhood $\Cal V_p$. By the result of Tian \cite{Ti}, we get the converging sequence of pseudo-metrics in the neighborhood $\Cal V_p$.  All the pseudo-metrics in the sequence are Bergman-type pseudo-metrics. Hence, we  get the convergence of the corresponding {\it  Hermitian kernel sections of positive type}\/  (\cite{A} and \cite {FK, definition on p.\,6}).  These authors consider Hermitian kernel functions. However, one can consider a positive Hermitian kernel section of the corresponding real bundle. Like the Bergman kernel of a domain, it is a reproducing kernel of the corresponding Hilbert space.

 As the limit, we obtain the diastasic potential  of the desired real analytic Kahler metric at the point $p$ because the corresponding $F$'s are holomorphic functional elements by \cite {FK, p.\, 6, Prop.\; I.1.1; p.\,12, Prop.\;I.1.6}. 

\example {$(3.4)$ Example}
Let $C$ be a nonsingular projective curve of genus $g(C)\geq 2$. We will assume: 
$$
\eusm K^{\ell}_C \subset \Cal L_C \subset \eusm K^{ m}_C
$$
where $\eusm K_C$ is the canonical bundle, and  $\ell, m$ are suitable  integers.
  We get the Bergman metrics on $C$ and Poincar\'e metric on $\Delta$.  It is well known \cite{Kr, Chap.\,III, Sect.\,4} 
that  the Bergman kernel  $\bold B_{\Delta,\eusm K^t}(z,\bar z)=c(t) \bold B^t_{\Delta,\eusm K}(z,\bar z)$, where $t\geq 1$ is an integer and $c(t)$ is a known constant. 
Hence 
$$
{1\over t} {\partial^2\log\bold B_{\Delta,\eusm K^t}(z,\bar z) \over \partial z \partial \bar z}= {\partial^2\log\bold B_{\Delta,\eusm K}(z,\bar z) \over \partial z \partial \bar z}.  
$$

Now, if $C$ is a general curvilinear section of  the projective variety  $X$ then we consider the inverse image of $C$ on $U_X$ with induced metric. We obtain an open Riemann surface $R_C \subset U_X$ in place of the disk $\Delta$.  It follows from the definition of $\Lambda$ that the restriction of $\Lambda$ on $R_C$ is its Bergman-Poincar\'e metric.
Indeed, high powers of $\Cal L_X |C$ are squeezed between powers of the canonical bundle on $C$.

Finally, the Bergman kernel of  
\endexample

\head
4.  Diastasic potentials generate a  function on $U_X$
\endhead

Let $z=(z_1, \dots, z_n)$ be a local coordinate system in a small neighborhood $\Cal V$ with origin at a fixed point $p_0\in \Cal V\subset U_X$. Let $\tilde \Phi_{p_0} (z(p), \overline{z(p)})$ be a diastasic potential at $p_0$. 

We consider
a path $[p_0, p]\subset U_X$, $p_0 =a_1, a_2, \dots, a_k=p$ where $a_i$ is close to $a_{i+1}$. Later we will impose additional restriction on the path. We would like to prolongate $\tilde \Phi_{p_0} (z(p), \overline{z(p)})$ along this path obtaining at each $a_i$ a diastasic potential $\tilde \Phi_{a_i} (z(p), \overline{z(p)})$ of our metric. 

We claim the prolongation is possible (and unique). Take the images of $a_i$ and $a_{i+1}$ in $X$.  Let  $C_{i,i+1}$ denote a curvilinear section of $X$ through the images.    We consider the inverse image of $C_{ij}$ on $U_X$ passing through the points $a_i$ and $a_{i+1}$, denoted by $R_{i,i+1}$. We will consider $R_{i,i+1}$ with the metric induced by $\Lambda$. We consider the geodesic on $R_{i,i+1}$ between the points $a_i$ and $a_{i+1}$. 

We obtain a broken geodesic (called also a broken/polygonal line) between the points $p_0$ and $p$. It remains to prolongate the $\tilde \Phi_{p_0} (z(p), \overline{z(p)})$ along this "broken geodesic" because these "broken geodesics" will approximate the original path.

%We get $g(C_{ij}) \geq 2$ since $\pi_1(X)$ is large. Then $R_{ij}:= \tau^{-1}
%(C_{ij})$ is a {\it connected}\/ open Riemann surface in $U_X$ by Deligne (see a %proof by Campana in \cite{Kol, Theorem 2.14.1}). 

We can assume the path $[a_i,a_{i+1}] \subset R_{ij}$. The restriction of $\Lambda$ on $R_{ij}$ is the Bergman-Poincar\'e metric whose diastasic potentials generate a function. Hence we can prolongate from $a_{i}$ to $a_{i+1}$. Here we use that the corresponding $F$ is a {\it holomorphic}\/ germ at a point on the diagonal in $U_X \times \bar{U}_X$ (Section 2). 

Since $\pi_1(U_X)$ is simply-connected we obtain the desired function $\bold P_U$ on $U_X$.

\head
5. End of proof of the theorem
\endhead
In view of the Oka-Grauert-Narasimhan theorem (Grauert's version), it remains to verify that for any real $\alpha$ the following set is relatively compact in $U_X$:
$$
\{q \in U_X \big |\bold P_U(q) < \alpha\}.
$$
Since $X$ is compact it is will suffice to consider a discrete set  $T$ in the fiber of $\tau$ over a point $Q\in X$.  We will show that $\bold P$ is unbounded on $T$. Take a general curvilinear section $C$ of $X$ through $Q$. Set
$
R_C:=\tau^{-1}(C) .
$ 
We obtain a connected open Riemann surface by Deligne's lemma (see a proof of the lemma by Campana in \cite{Kol, Theorem 2.14.1}). In particular, the natural map $\pi_1(R_C) \rightarrow \pi_1(X)$ is surjective hence $g(C)\geq 2$.

By the fundamental property of the diastasis (Sect.\;2), $\bold P_R = \bold P _U|R_C$ where $\bold P_R$ is the corresponding function generated by diastasic  potentials on $R_C$. The restriction of $\Lambda$ on $R_C$ is its Bergman-Poincar\'e metric. Hence 
the corresponding set
$$\{q \in R_C \; |\; \tau(q)=Q,\; \bold P_R(q) < \alpha\}
$$
 is finite like in Siegel's proof (see Sect.\;1 and \cite{Kob1, Theorem 9.5}).

This proves the theorem.

The metric was suggested by a problem of Yau \cite{Y} (also, see \cite {Ti}), and it is a generalization of the classical Poincar\'e metric on the disk $\Delta$. The construction in (3.1) is explained in the paper by Tian \cite{Ti}.

(3.1)
Let $C$ be a nonsingular projective curve of genus $g(C)\geq 2$. We will assume: 
$$
\eusm K^{\ell}_C \subset \Cal L_C \subset \eusm K^{ m}_C
$$
where $\eusm K_C$ is the canonical bundle, and  $\ell, m$ are suitable  integers.
  We get the Bergman-Poincar\'e metrics on $C$ and the Poincar\'e metric on $\Delta$.  It is well known \cite{Kr, Chap.\,III, Sect.\,4} 
that  the Bergman kernel  $\bold B_{\Delta,\eusm K^t}(z,\bar z)=c(t) \bold B^t_{\Delta,\eusm K}(z,\bar z)$, where $t\geq 1$ is an integer and $c(t)$ is a known constant depending on $t$ only. 
Hence 
$$
{1\over t} {\partial^2\log\bold B_{\Delta,\eusm K^t}(z,\bar z) \over \partial z \partial \bar z}= {\partial^2\log\bold B_{\Delta,\eusm K}(z,\bar z) \over \partial z \partial \bar z}.  
$$

Now, if $C$ is a general curvilinear section of  the projective variety  $X$ then we consider the inverse image of $C$ on $U_X$ with induced metric. We obtain an open Riemann surface $R_C \subset U_X$ in place of the disk $\Delta$.  It follows from the definition of $\Lambda$ that the restriction of $\Lambda$ on $R_C$ is its Bergman-Poincar\'e metric.
Indeed, high powers of $\Cal L_X |C$ are squeezed between powers of the canonical bundle on $C$.

\smallskip
{\it  Acknowledgment.} The author would like to thank Pierre Milman for very helpful conversations, and the referee for very helpful remarks.

\enddocument

\ref  
\key D  \by  D. Drasin \pages 356--365
\paper Cusp forms and Poincar\'e series
\yr1968 \vol  90
\jour Amer. J.  Math.
\endref

\ref  
\key  U  
\by M. Umehara \pages  203--214
\paper  Kaehler Submanifolds of Complex Space Forms
\jour Tokyo J. Math
\vol 10, No. 1
\endref
\endRef

\enddocument